\documentclass[12pt,a4paper]{amsart}

\usepackage[left=25truemm,right=25truemm,top=30truemm,bottom=30truemm]{geometry}

\usepackage{mathtools}
\usepackage{amssymb}
\usepackage{thmtools}
\usepackage{mathrsfs}

\usepackage[dvipdfmx]{graphicx}
\usepackage{subcaption}
\usepackage{float}
\usepackage[dvipsnames]{xcolor}
\usepackage{framed}
\usepackage[all]{xy}
\usepackage{listings}
\usepackage{setspace}
\usepackage[english]{babel}
\usepackage{amscd}

\usepackage{tikz}
\usetikzlibrary{trees}

\usepackage{hyperref}
\usepackage{cleveref}

\theoremstyle{plain}
\newtheorem{thm}{\textrm{Theorem}}[section]
\newtheorem{lem}[thm]{\textrm{Lemma}}
\newtheorem{cor}[thm]{\textrm{Corollary}}
\newtheorem{prop}[thm]{\textrm{Proposition}}

\theoremstyle{definition}

\newtheorem{defi}[thm]{\textrm{Definition}}
\newtheorem{rem}[thm]{\textrm{Remark}}

\DeclareMathOperator{\Diff}{\mathop{\rm{Diff}}}
\DeclareMathOperator{\Int}{\mathop{\rm{Int}}}

\Crefname{figure}{Figure}{Figures}

\numberwithin{equation}{section}

\newcommand{\id}{\rm{id}}

\title{\small Topology of boundary special generic maps\\ into Euclidean spaces}

\author[K. Iwakura]{Koki Iwakura}
\address{Joint Graduate School of Mathematics for Innovation, Kyushu University, Motooka 744, Nishiku, Fukuoka 819-0395, Japan.}
\email{iwakura.koki0105@gmail.com}

\date{}

\begin{document}

\maketitle

\begin{abstract}
We introduce boundary special generic maps, a class of submersions from manifolds with boundary to Euclidean spaces whose restriction to the boundary has only boundary definite fold points as its singular points. 
We derive the differential-topological restrictions imposed by the existence of such maps on the global structure of the source manifolds. 
Furthermore, we apply our results to the non-singular extension problem, which asks when a map on a closed manifold extends to a non-singular map on a manifold with boundary, and obtain new results on non-singular extensions of special generic maps.
\end{abstract}

\section{Introduction}

\subsection{Main results on boundary special generic maps}
Singular points of a smooth map are local objects; however, numerous works have shown that they can impose constraints on the global structure of the source manifold, including its topology and smooth structure (for example, see~\cite{Sae1, SS1, SS2, Tho}). 
Many previous works have focused on maps defined on manifolds without boundary, aiming to understand how singularities reflect or restrict global structure. 
In contrast, for maps whose source manifolds have non-empty boundary, such a framework is far less developed: while the case of functions is well understood through Morse theory for manifolds with boundary (for details, see~\cite{BNR, Bra, Haj, JR}), comparable results for general maps remain limited. 
Indeed, in the presence of boundary, one must control not only singular points in the interior of the manifold but also those arising from the restriction of the map to the boundary. 
Consequently, the methods available for closed manifolds do not carry over in any straightforward way to the case of manifolds with boundary.

\vspace{2pt}

We accordingly introduce boundary special generic maps, a class of smooth maps from compact, connected manifolds with non-empty boundary to Euclidean spaces that have no singular points in the interior and whose restrictions to the boundary admit only boundary definite fold points as singular points. 
Our aim is to derive differential-topological restrictions on the source manifolds from the properties of these singular points. 
This notion was introduced by Shibata~\cite{Shi} in the case of maps from $3$-manifolds with boundary into the plane, and here we formulate a general definition in higher dimensions. 
Moreover, boundary special generic maps can be regarded as a natural generalization of the following Morse-theoretic situation, which occurs frequently as the simplest case: a smooth function on a manifold with boundary has no critical points in the interior, while its restriction to the boundary has critical points; each such point is a local maximum or minimum of the original function.
Our main results are summarized in the following four theorems.

\vspace{2pt}

\begin{thm}%定理1.1
Let $N$ be a compact, connected $n$-dimensional manifold with boundary, where $n\geq 2$. 
Then, $N$ admits a boundary special generic map into $\mathbb{R}$ if and only if $N$ is diffeomorphic to $D^n$. 
\end{thm}

\vspace{0.05\baselineskip}

Theorem\,1.1 determines the diffeomorphism type of the source manifold of a boundary special generic map into $\mathbb{R}$. 
This result may be regarded as a natural analogue of Reeb's sphere theorem~\cite{Mil1} for manifolds with boundary. 
The proof relies on the classical results on the diffeomorphism group of the disk~\cite{Hat, Sma} and standard consequences of the h-cobordism theorem \cite{Mil1}.

\vspace{2pt}

By contrast, when the target is $\mathbb{R}^2$, the above arguments do not apply in any straightforward way. 
We instead use properties of the Reeb space, obtained by collapsing each connected component of a fiber to a point, to prove the following result.

\vspace{2pt}

\begin{thm}%定理1.2
Let $N$ be a compact, connected $n$-dimensional manifold with boundary, where $n\geq 3$. 
Then, $N$ admits a boundary special generic map into $\mathbb{R}^2$ if and only if $N$ is diffeomorphic to a boundary sum of finitely many $D^{n-1}$-bundles over $S^1$.
Here, we interpret the case of zero summands as $D^n$.
\end{thm}

\vspace{0.05\baselineskip}

Theorem\,1.2 determines the diffeomorphism type of the source manifold of a boundary special generic map into $\mathbb{R}^2$. 
The key point of the proof is a handle decomposition of the source manifold induced from a handle decomposition of the Reeb space.

\vspace{2pt}

Moreover, the method used in the proof of Theorem\,1.2 also applies when the target manifold is $\mathbb{R}^3$, yielding the following result.

\vspace{2pt}

\begin{thm}%定理1.3
Let $N$ be a compact, simply connected $n$-dimensional manifold with boundary, where $n\geq 4$ and $n\neq 6,7$. 
If $N$ admits a boundary special generic map into $\mathbb{R}^3$, then $N$ is diffeomorphic to a boundary sum of finitely many $D^{n-2}$-bundles over $S^2$. 
Here, we interpret the case of zero summands as $D^n$.
\end{thm}

\vspace{0.05\baselineskip}

Theorem\,1.3 describes the possible diffeomorphism type of the source manifold of a boundary special generic map into $\mathbb{R}^3$, except in several dimensions. 
The excluded dimensions require separate arguments and will be treated individually.

\vspace{2pt}

As the dimension of the target increases, the handle decomposition of the Reeb space becomes increasingly complicated, and the methods used for Theorems\,1.2 and 1.3 do not extend in any straightforward way to targets of higher dimension. 
Therefore, we introduce an alternative approach that is effective across many dimensions: we decompose the source manifold of a boundary special generic map into two disk bundles, corresponding to the neighborhood of the set of boundary definite fold points and the complement of this neighborhood. 
This decomposition can be viewed as an analogue of the handle decompositions for Morse functions on manifolds with boundary developed in~\cite{BNR, JR}.

\vspace{2pt}

\begin{thm}%定理1.4
Let $N$ be a compact, connected $n$-dimensional manifold with boundary, where $n\geq 6$.
Assume that $n-k=1$. 
Then, $N$ admits a boundary special generic map into $\mathbb{R}^k$ with contractible Reeb space if and only if $N$ is diffeomorphic to $D^n$. 
\end{thm}

\vspace{0.05\baselineskip}

Theorem\,1.4 determines the diffeomorphism type of the source manifold of a boundary special generic map when the target has dimension at least five, and the codimension is small, under the additional assumption that the Reeb space is contractible. 
The proof relies on classical results on the diffeomorphism group of the disk~\cite{Hat, Sma}.
Taken together, Theorems 1.1--1.4 yield restrictions on the diffeomorphism types of source manifolds of boundary special generic maps across various dimensions of the targets. 
We now apply these results to the non-singular extension problem.

\vspace{2pt}

\subsection{Applications to the non-singular extension problem}
The problem of finding necessary and sufficient conditions for a given map on a closed manifold to admit an extension to a non-singular map on a manifold with the given manifold as boundary is called the non-singular extension problem (see~\cite{Cur, Iwa, Sei}). 
It can be viewed as an analogue, in the case where the dimension of the source manifold exceeds that of the target manifold, of the classical extension problems for immersions (see~\cite{Poe, Mar}), and has been studied in the context of differential topology and singularity theory. 
On the other hand, controlling singular behavior that may occur throughout the manifold and ruling out the appearance of singular points along an extension is highly challenging; consequently, for maps between manifolds in general dimensions, a systematic understanding remains limited.

\vspace{2pt}

In this paper, we focus on the non-singular extension problem for special generic maps on closed manifolds. 
A special generic map is a smooth map whose singular points are all definite fold points, the natural higher-dimensional generalization of local maxima and minima. 
Nevertheless, even for such basic maps, systematic results on non-singular extensions are scarce. 
Since the restriction to the boundary of a boundary special generic map is a special generic map, we can exploit the results of Section\,1.1 to obtain the following four consequences for non-singular extensions of special generic maps.

\vspace{2pt}

\begin{cor}
Let $M$ be a closed $m$-dimensional manifold, where $m\geq 1$, and $f\colon M\to\mathbb{R}$ be a special generic map. 
Then, $f$ admits a non-singular extension if and only if $M$ is diffeomorphic to $S^m$. 
\end{cor}

\vspace{0.05\baselineskip}

Corollary\,1.5 gives necessary and sufficient conditions for a special generic map into $\mathbb{R}$ to admit a non-singular extension. 
The proof combines Theorem\,1.1 with a result of Seigneur~\cite{Sei} on non-singular extensions of Morse functions on the sphere. 
In particular, it follows that any special generic map from an exotic sphere to $\mathbb{R}$ does not admit a non-singular extension. 
Motivated by this, we next raise the dimension of the target and investigate similar phenomena, namely that special generic maps on certain manifolds admit no non-singular extension by boundary special generic maps.

\vspace{2pt}

\begin{cor}
Special generic maps from $\Sigma^m$ or $(S^1\times S^{m-1})\,\sharp\,\Sigma^m$ to $\mathbb{R}^2$ do not admit boundary special generic maps as non-singular extensions, where $\Sigma^m$ is an $m$-dimensional exotic sphere with $m\geq 7$. 
\end{cor}

\vspace{0.05\baselineskip}

It is known that an exotic sphere $\Sigma^m$, as well as its connected sums with $S^1\times S^{m-1}$, admits a special generic map into $\mathbb{R}^2$ by~\cite{Sae1}.
However, the existence of the above non-singular extension imposes restrictions on the smooth structure of the source manifold; in particular, it excludes $\Sigma^m$ and $(S^1\times S^{m-1})\,\sharp\,\Sigma^m$ as possible source manifolds of a special generic map admitting such an extension.
The proof combines Theorem\,1.2 with Schultz's computation~\cite{Sch} of the inertia group of $S^1\times S^{m-1}$.

\vspace{2pt}

We then turn to special generic maps into $\mathbb{R}^3$. 
In this case, we restrict our attention to $3$-manifolds as the source manifolds, where specific features of $3$-manifolds can be utilized.

\vspace{2pt}

\begin{cor}
Special generic maps from a closed, connected, irreducible $3$-manifold with non-cyclic fundamental group to $\mathbb{R}^3$ admit no boundary special generic map as a non-singular extension whose source manifold is simply connected. 
\end{cor}

\vspace{0.05\baselineskip}

It is known that every closed, orientable $3$-manifold admits a special generic map into $\mathbb{R}^3$ by \cite{Eli}.
However, the existence of such a non-singular extension imposes strong restrictions on the fundamental group of the $3$-manifold; in particular, it rules out irreducible $3$-manifolds with non-cyclic fundamental group as possible source manifolds of a special generic map admitting such an extension.
The proof uses Theorem\,1.3 together with basic properties of lens spaces.
Moreover, using this corollary, we obtain Corollary\,3.14, showing that special generic maps on Brieskorn manifolds into $\mathbb{R}^3$ admit no non-singular extensions to Mazur-type manifolds by boundary special generic maps.

\vspace{2pt}

We next consider the case of target manifolds with general dimensions. 
Using the fact by~\cite{KM} that exotic spheres are stably parallelizable, we obtain the following consequence.

\vspace{2pt}

\begin{cor}
Special generic maps from $\Sigma^m$ to $\mathbb{R}^m$ do not admit boundary special generic maps with the contractible Reeb space as non-singular extensions, where $\Sigma^m$ is an $m$-dimensional exotic sphere with $m\geq 7$. 
\end{cor}

\vspace{0.05\baselineskip}

It is known that an exotic sphere admits a special generic map into the Euclidean space $\mathbb{R}^m$ of the same dimension~\cite{Eli}, since it is stably parallelizable. 
However, the existence of such a non-singular extension imposes strong restrictions on the smooth structure of the source manifold; in particular, it excludes $\Sigma^m$ as a possible source manifold of a special generic map admitting such an extension.
The proof follows directly from Theorem\,1.4.

\vspace{2pt}

Taken together, these results indicate that the existence of non-singular extensions can impose stronger constraints on closed manifolds than the mere existence of special generic maps. 
In particular, Corollaries\,1.5, 1.6, and 1.8 relate non-extendibility to smooth structures of closed manifolds, while Corollary\,1.7 yields a restriction on the fundamental group in dimension three. 
As a consequence of these constraints, many previously known examples of special generic maps do not admit a boundary special generic map as a non-singular extension. 
Moreover, to our knowledge, a systematic framework for producing explicit examples without extendibility has not appeared in the existing literature.

\vspace{2pt}

\subsection{Contents}
Section\,2 is devoted to the main results on boundary special generic maps (Theorems\,1.1--1.4). 
In Subsection\,2.1, we recall basic notions about singular points of smooth maps and introduce the Reeb space of maps. 
In Subsection\,2.2, we introduce boundary special generic maps. 
In Subsections\,2.3--2.5, we prove Theorems\,1.1--1.3, corresponding to the cases where the target manifolds are $\mathbb{R}$, $\mathbb{R}^2$, and $\mathbb{R}^3$, respectively. 
Subsection\,2.6 deals with the case of the target manifold with general dimensions and proves Theorem\,1.4. 
In Section\,3, we formulate the non-singular extension problem and prove the application results (Corollaries\,1.5--1.8). 
Subsections\,3.1--3.2 introduce the notion of non-singular extensions and special generic maps. 
Then, Subsections\,3.3--3.6 prove Corollaries\,1.5--1.8 according to the dimensions of the target manifolds.

\vspace{2pt}

\subsection{Notations}
Throughout this paper, unless otherwise stated, all manifolds and maps between them are assumed to be smooth of class $C^{\infty}$. 
We use the symbol $\cong$ to denote diffeomorphism between manifolds and isomorphism between groups or rings. 
Furthermore, $\Diff$ denotes the diffeomorphism group of the relevant manifold; when the manifold has non-empty boundary, we do not impose conditions on the behavior along the boundary.

\vspace{2pt}

\section{On boundary special generic maps}
In this section, we aim to prove Theorems\,1.1--1.4. 
We begin with a recall of basic terminology from singularity theory that will be used throughout the paper.
We then introduce boundary special generic maps and give some standard facts about their Reeb spaces for later arguments. 
Finally, we prove Theorems\,1.1--1.4, beginning with Theorem\,1.1.

\vspace{2pt}

\subsection{Basic terminology from singularity theory}
In this subsection, we briefly review the basic notations used in the paper; see \cite{GG} and \cite{Sae2} for further details.
Throughout this subsection, let $M$ and $K$ be manifolds (possibly with non-empty boundary) with $\dim M\geq \dim K\geq 1$, and let $f \colon M \to K$ be a map.

\vspace{2pt}

\begin{defi}[\textbf{Singular point}]
A point $p \in M$ is called a \textit{singular point} of $f$ if $\operatorname{rank} df_p < \dim K$. 
Otherwise, $p$ is called a \textit{regular point} of $f$. 
In particular, $f$ is called a \textit{submersion} if every point of $M$ is a regular point of $f$.
\end{defi}

\vspace{2pt}

We denote by $S(f)$ the subset of $M$ consisting of all singular points of $f$:
$$
S(f) = \{\, p \in M \mid \operatorname{rank} df_p < \dim K \,\}.
$$
We call $S(f)$ the \textit{singular point set} of $f$.

\vspace{2pt}

We next recall the definition of the Reeb space of $f$, obtained by collapsing each connected component of each fiber to a point.
This space will be used repeatedly later, as it provides a convenient way to encode the global structure of a map.

\vspace{2pt}

\begin{defi}[\textbf{Reeb space}]
For $x, y\in M$, we define a relation $``\sim"$ on $M$ by
$$
x \sim y \iff  f(x) = f(y) \text{ and $x,y$ lie in the same connected component of } f^{-1}(f(x)).
$$
Since this relation is an equivalence relation on $M$, the quotient space $W_f=M/\sim$ is a topological space equipped with the quotient topology, and is called the \textit{Reeb space} of $f$. 
Moreover, there exists a unique continuous map $\overline{f}$, called the \textit{Reeb map} of $f$, that makes the following diagram commutative, 
$$
       \xymatrix{
         M \ar[r]^{f} \ar[d]_{q_{f}} & K \\
         W_{f} \ar[ru]_{\overline{f}},
       }
$$
where, $q_f\colon M\to W_f$ denotes the quotient map. 
\end{defi}

\vspace{2pt}

\subsection{Definition of boundary special generic maps} 
In this subsection, we introduce boundary special generic maps. 
In the rest of this subsection, let $N$ be a compact, connected $n$-dimensional manifold with boundary, $K$ be a $k$-dimensional open manifold with $n>k$, and $F\colon N\to K$ be a map.

\vspace{2pt}

\begin{defi}[\textbf{Boundary definite fold point}]
A point $p\in\partial N$ is called a \textit{boundary definite fold point} of $F$ if there exist local coordinates $(x_1,\dots, x_n)$ of $N$ around $p$ and $(y_1,\dots, y_k)$ of $K$ around $F(p)$ such that $F$ has the local form 
\begin{equation*}
F; (x_1,\cdots,x_n)\mapsto(y_1,\cdots,y_k)
=\Big(x_1,\cdots,x_{k-1},\sum_{i=k}^{n-1} x_i^2+x_n\Big), 
\end{equation*} 
where $x_n>0$ and $x_n=0$ correspond to $\Int N$ and $\partial N$, respectively. 
\end{defi}

\begin{defi}[\textbf{Boundary special generic map}]
A \textit{boundary special generic map} is a submersion $F\colon N\to K$ such that the singular points of $F|_{\partial N}$ are all boundary definite fold points of $F$. 
\end{defi}

\begin{rem}
In the above definitions, the boundary definite fold points of $F$ are not singular points of $F$ itself.  
However, since $x_n=0$ corresponds to $\partial N$, these points are the singular points of $F|_{\partial N}$.
In particular, they are definite fold points of $F|_{\partial N}$, and hence $F|_{\partial N}$ is a special generic map; see Definitions\,3.3 and~3.4.
\end{rem}

\vspace{2pt}

In this paper, we consider only boundary special generic maps into Euclidean spaces. 
Namely, throughout the paper, we deal with maps of the form $F \colon N \to \mathbb{R}^k$.

\vspace{2pt}

Although we will not use the next statement immediately, we record it here for later reference.
Since the proof is analogous to that in \cite{Sae2}, we omit it.

\vspace{2pt}

\begin{prop}
The Reeb space $W_F$ is a compact, connected, orientable $k$-dimensional manifold with boundary. 
Moreover, the Reeb map $\overline{F}\colon W_F\to\mathbb{R}^k$ is an immersion, and the quotient map $q_F\colon N\to W_F$ is a smooth map, in such a way that the restriction $q_F|_{S(F|_{\partial N})}\colon S(F|_{\partial N})\to \partial W_F$ to the singular point set $S(F|_{\partial N})$ is a diffeomorphism. 
\end{prop}

\vspace{2pt}

\subsection{Proof of Theorem\,1.1}
In this subsection, we prove Theorem\,1.1.
Assume that there exists a boundary special generic map $F\colon N\to \mathbb{R}$.
By the definition of boundary definite fold points, the manifold $N$ is obtained by gluing two copies of $D^n$ along a diffeomorphism of $D^{n-1}$ on their boundaries. 
In particular, $N$ is homeomorphic to $D^n$. 
We show that such a manifold $N$ is diffeomorphic to $D^n$ for every $n\geq 2$.

\vspace{2pt}

We first consider the cases $n=2,3$.
Since $N$ is homeomorphic to $D^n$, it is diffeomorphic to $D^n$ by the uniqueness of the smooth structures in dimensions $2$ and~$3$.

\vspace{2pt}

We next consider the case $n=4$. 
By~\cite{Hat}, we have a homotopy equivalence $\Diff(D^3)\simeq {\rm{O}}(3)$.
In particular, $\Diff(D^3)$ has two components, represented by the identity and a reflection. 
Hence, the gluing diffeomorphism of $D^{3}$ is isotopic to either the identity or a reflection.
Therefore, $N$ is diffeomorphic to $D^{4}$.

\vspace{2pt}

We next consider the case $n=5$.
The restriction $F|_{\partial N}$ has exactly two critical points, a maximum and a minimum.
By Reeb's sphere theorem~\cite{Mil1}, $\partial N$ is obtained by gluing two copies of $4$-dimensional disks along their boundaries.
In particular, $\partial N$ is diffeomorphic to $S^4$ by~\cite{Cer}.
Since $N$ is homeomorphic to $D^5$, it follows from~\cite[Proposition C, p.110]{Mil2} that $N$ is diffeomorphic to $D^5$.

\vspace{2pt}

Finally, we consider the case $n\geq 6$. 
By Reeb's sphere theorem~\cite{Mil1}, $\partial N$ is homeomorphic to $S^{n-1}$.
In particular, $\partial N$ is simply connected. 
Since $N$ is homeomorphic to $D^n$, it follows from~\cite[Proposition A,~p.108]{Mil2} that $N$ is diffeomorphic to $D^n$.

\vspace{2pt}

We next prove the converse direction of Theorem\,1.1 by constructing a boundary special generic map from $N$ to $\mathbb{R}$. 
Since $N$ is diffeomorphic to $D^n$, such a map is obtained as the composition
$$
D^{n}\hookrightarrow\mathbb{R}^{n}\xrightarrow{\pi}\mathbb{R}, 
$$
where the first map is the natural inclusion and $\pi$ is the standard projection. 
This completes the proof of Theorem\,1.1.
\qed

\vspace{2pt}

\subsection{Proof of Theorem\,1.2}
In this subsection, we prove Theorem\,1.2. 
Assume first that there exists a boundary special generic map $F\colon N\to\mathbb{R}^2$.
We show that a handle decomposition of $W_F$ induces, via $q_F$, a handle decomposition of $N$ that consists only of a $0$-handle and $1$-handles, which implies that $N$ has the required diffeomorphism type.

\vspace{2pt}

By Proposition\,2.6, the Reeb space $W_F$ is a compact, connected, orientable surface with boundary.
Fix a handle decomposition
\begin{equation}
W_F=h^0\cup h^1_1\cup\dots\cup h^1_R. 
\end{equation}
Here, $h^0$ is a $0$-handle and $h^{1}_i$ is a $1$-handle for $i=1,\dots,R$.
The integer $R\geq 0$ is the number of $1$-handles. 
If $R=0$, we regard $W_F$ as $h^0$. 
For each $i$, let $I_i$ and $I'_i$ denote the attaching regions of $h^1_i$. 
Note that these regions are pairwise disjoint in $\partial h^0$.

\vspace{2pt}

We first consider the case $R=0$. 
Then, $W_F=h^0$ and $h^0\cong D^2$. 
Let $h\colon D^2\to\mathbb{R}$ be the height function.
Then, the composition
$$
N=q_F^{-1}(h^0)\xrightarrow{q_F} h^0\cong D^2\xrightarrow{h}\mathbb{R}
$$ 
is a boundary special generic map.
Therefore, $N$ is diffeomorphic to $D^n$ by Theorem\,1.1.

\vspace{2pt}

Next, assume that $R\geq 1$. 
We claim that, for each $i=1,\dots,R$,  
$$
q_F^{-1}(I_i)\cong D^{n-1}, 
\qquad
q_F^{-1}(I'_i)\cong D^{n-1}, 
$$
and
$$
q_F^{-1}(h^0)\cong D^{n}, 
\qquad
q_F^{-1}(h^1_i)\cong D^{1}\times D^{n-1}. 
$$
We briefly explain these identifications. 
First, for each attaching region $I_i\subset \partial h^0$, the composition
$$
q_F^{-1}(I_i)\xrightarrow{q_F} I_i\hookrightarrow \mathbb{R}
$$
is a boundary special generic map, where the second map is an embedding.
Hence, $q_F^{-1}(I_i)\cong D^{n-1}$ by Theorem\,1.1.
The same argument gives $q_F^{-1}(I_i')\cong D^{n-1}$.
Next, the manifold $q_F^{-1}(h^0)$ is obtained by gluing copies of $D^1\times D^{n-1}$ to $D^2\times D^{n-2}$ along ${\id}_{D^{n-1}}$, corresponding to the components of $\overline{\partial h^0\setminus \sqcup_{i=1}^{R}{I_i\cup I'_{i}}}$. 
Thus, $q_F^{-1}(h^0)\cong D^{n}$. 
Finally, since $q_F^{-1}(h^1_i)$ is a $q_F^{-1}(I_i)$-bundle over $D^1$ and $q_F^{-1}(I_i)\cong D^{n-1}$, we have $q_F^{-1}(h^1_i)\cong D^{1}\times D^{n-1}$.

\vspace{2pt}

By~(2.1) and the above identifications, $N$ admits the handle decomposition induced by $q_F$.
In particular, $N$ is an $n$-dimensional $1$-handlebody.
Hence, $N$ is obtained as a boundary sum of $R$ $D^{n-1}$-bundles over $S^1$.

\vspace{2pt}

Conversely, assume that $N$ is a boundary sum of finitely many $D^{n-1}$-bundles over $S^1$.
Then, there exist integers $r,r'\geq 0$ such that 
$$
N\cong(\natural^r S^1\times D^{n-1})\,\natural\,(\natural^{r'} S^{1}\tilde{\times} D^{n-1}), 
$$
where $S^{1}\tilde{\times} D^{n-1}$ is the non-trivial one. 
We construct a boundary special generic map of this manifold into $\mathbb{R}^2$ as follows.

\vspace{2pt}

We first consider the case $r=r'=0$. 
Then, $N\cong D^n$. 
In this case, the composition
$$
D^n \hookrightarrow \mathbb{R}^n \xrightarrow{\pi} \mathbb{R}^2
$$
is a boundary special generic map, where the first map is the natural inclusion and $\pi$ is the standard projection.

\vspace{2pt}

Assume next that $r+r'\geq 1$. 
We use the following lemma, which can be verified by an argument similar to that of~\cite[Lemma 5.4]{Sae1}.

\vspace{2pt}

\begin{lem}
Let $N_j$ be a compact, connected $n$-dimensional manifold with boundary that admits a boundary special generic map into $\mathbb{R}^2$ for $j=1,2$. 
Then, there exists a boundary special generic map from $N_1 \,\natural\, N_2$ to $\mathbb{R}^2$.
\end{lem}

\vspace{2pt}

By this lemma, it suffices to construct boundary special generic maps into $\mathbb{R}^2$ for the two building blocks $S^1\times D^{n-1}$ and $S^1\tilde{\times}D^{n-1}$.

\vspace{2pt}

A boundary special generic map from $S^1\times D^{n-1}$ to $\mathbb{R}^2$ is given by
$$
S^1 \times D^{n-1} 
\xrightarrow{\,{\id}_{S^1} \times h\,} 
S^1 \times \mathbb{R} 
\hookrightarrow 
\mathbb{R}^2, 
$$
where $h\colon D^{n-1}\to\mathbb{R}$ is the height function and the second map is an embedding.

\vspace{2pt}

Similarly, a boundary special generic map from $S^1\tilde{\times}D^{n-1}$ to $\mathbb{R}^2$ is given by  
$$
S^1 \,\tilde{\times}\, D^{n-1} 
\xrightarrow{\,\tilde{h}\,} 
S^1 \times \mathbb{R} 
\hookrightarrow 
\mathbb{R}^2, 
$$
where $\tilde{h}([t,x_1,\dots,x_{n-1}])=(t,x_{n-1})$ with the identification
$$
S^1\tilde{\times}D^{n-1}=[-1,1]\times D^{n-1}/(1,(x_1,\cdots,x_{n-2},x_{n-1}))\sim (-1,(x_1,\cdots,x_{n-2},-x_{n-1})), 
$$
and, the second map is an embedding.

\vspace{2pt}

Applying the lemma repeatedly, we obtain a boundary special generic map from $N$ to $\mathbb{R}^2$.
This completes the proof of Theorem\,1.2.
\qed

\vspace{2pt}

\subsection{Proof of Theorem\,1.3}
In this subsection, we prove Theorem\,1.3. 
Let $F\colon N\to\mathbb{R}^3$ be a boundary special generic map.
As in the proof of Theorem\,1.2, we show that a handle decomposition of $N$ is induced from a handle decomposition of $W_F$ via $q_F$. 
This yields a handle decomposition of $N$ consisting only of a $0$-handle and $2$-handles, which implies the claimed diffeomorphism type in Theorem\,1.3.

\vspace{2pt}

Since $N$ is simply connected, the Reeb space $W_F$ is simply connected (see Remark~2.15 in subsection~2.6). 
Hence, $W_F$ is diffeomorphic to the $3$-manifold obtained from $D^3$ by removing finitely many open $3$-dimensional balls.
We denote the number of such balls by $R\geq 0$. 
Then, $W_F$ admits a handle decomposition 
\begin{equation}
W_F=h^0\cup h^2_1\cup\cdots\cup h^2_R,
\end{equation}
where $h^0$ is a $0$-handle and each $h^2_i$ is a $2$-handle.
For each $i$, let $A_i\subset \partial h^0$ be the attaching region of $h^2_i$.
Note that the regions $A_1, \dots, A_R$ are pairwise disjoint in $\partial h^0$.

\vspace{2pt}

First, we consider the case $R=0$. 
Then, $W_F=h^0$ and $h^0\cong D^3$. 
Let $h\colon D^3\to\mathbb{R}$ be the height function.
Then, the composition 
$$
N=q_F^{-1}(h^0)\xrightarrow{q_F} h^0\cong D^3\xrightarrow{h}\mathbb{R}, 
$$ 
is a boundary special generic map. 
Therefore, $N$ is diffeomorphic to $D^n$ by \textrm{Theorem\,1.1}.

\vspace{2pt}

Next, assume that $R\geq 1$. 
We claim that, for each $i=1,\dots,R$, 
$$
q_F^{-1}(h^2_i)\cong D^{2}\times D^{n-2}, 
\qquad
q_F^{-1}(A_i)\cong S^1\times D^{n-2}, 
\qquad
q_F^{-1}(h^0)\cong D^{n}. 
$$
We briefly explain these identifications.
First, by Theorem\,1.1, the inverse image of the cocore of $h^2_i$ is diffeomorphic to $D^{n-2}$.
Hence, $q_F^{-1}(h^2_i)$ is a $D^{n-2}$-bundle over $D^2$, and therefore
$q_F^{-1}(h^2_i)\cong D^2\times D^{n-2}$.
In particular, restricting to the attaching region of $h^2_i$ gives $q_F^{-1}(A_i)\cong S^1\times D^{n-2}$.
Finally, the manifold $q_F^{-1}(h^0)$ is obtained from $D^3\times D^{n-3}$ by attaching a $D^{n-2}$-bundle over each component $S$ of $\overline{\partial h^0\setminus(\sqcup_{i=1}^{R} A_i)}$ along the identity on $S\times D^{n-3}$. 
Each such $S$ is a compact, connected, orientable surface of genus $0$.
Hence, every $D^{n-2}$-bundle over $S$ is trivial.
(Here, the assumption $n\neq 6,7$ is used; see Remark\,2.8.)
Therefore, $q_F^{-1}(h^0)$ is diffeomorphic to $D^n$.

\vspace{2pt}

By~(2.2) and the above identifications, $N$ admits the handle decomposition induced by $q_F$.
In particular, $N$ is an $n$-dimensional $2$-handlebody.
Since $N$ is orientable, the framing of each $2$-handle corresponds to an element of $\pi_1({\rm{SO}}(n-2))$. 
When $n\geq 5$, $\pi_1({\rm{SO}}(n-2))\cong\mathbb{Z}_2$; when $n=4$, $\pi_1({\rm{SO}}(2))\cong\mathbb{Z}$. 
Therefore, if $n\geq 5$, $N$ is a boundary sum of $R$ $D^{n-2}$-bundles over $S^{2}$ by a standard fact.
If $n=4$, we also need to consider the positions of the attaching circles in $\partial h^0\cong S^3$.
In this case, by pulling back properly embedded disks in $W_F$ that bound the attaching circles of the $2$-handles via $q_F$, we see that these attaching circles form an $R$-component unlink in $S^3$.
Hence, $N$ is a boundary sum of $R$ $D^{2}$-bundles over $S^{2}$.
\qed

\vspace{2pt}

\begin{rem}
In the proof above, the assumption $n\neq 6,7$ is used to show that the $D^{n-2}$-bundle over $S$ is trivial. 
If $\sharp S(F|_{\partial N})\leq 2$, then each component $S$ of $\overline{\partial h^0\setminus(\sqcup_{i=1}^{R} A_i)}$ is diffeomorphic to $D^2$. 
Hence, the similar conclusion also holds for $n=6,7$ under the assumption $\sharp S(F|_{\partial N})\leq 2$. 
Here, $\sharp$ denotes the number of connected components of the relevant set. 
\end{rem}

\begin{rem}
For the statement of Theorem\,1.3, when $n\ge 5$, the converse direction also holds.
That is, for a boundary sum of $D^{n-2}$-bundles over $S^2$, there exists a boundary special generic map into $\mathbb{R}^3$.
This can be proved by an argument similar to that in the proof of Theorem\,1.2.
On the other hand, it is not known whether every boundary sum of $D^2$-bundles over $S^2$ admits a boundary special generic map into $\mathbb{R}^3$.
\end{rem}

\vspace{2pt}

\subsection{Proof of Theorem\,1.4}
In this subsection, we prove Theorem\,1.4.
We first state the following lemma, which decomposes the source manifold of a boundary special generic map into disk bundles.
Since it can be proved by essentially the same argument as in~\cite{Sae1}, we omit the proof.

\vspace{2pt}

\begin{lem}\label{lem:collar}
Let $N$ be a compact, connected $n$-dimensional manifold with boundary, and let $F\colon N\to\mathbb{R}^k$ be a boundary special generic map.
Let $C$ be a collar neighborhood of $\partial W_F\subset W_F$.
Identify $C$ with $\partial W_F\times[0,1]$ so that $\partial W_F$ corresponds to $\partial W_F\times\{0\}$.
Put $W=\overline{W_F\setminus C}$.
Then:
\begin{enumerate}
\setlength{\leftskip}{-14pt}
\item $q_F|_{q_F^{-1}(W)} \colon q_F^{-1}(W) \to W$ is a $D^{n-k}$-bundle with structure group $\Diff(D^{n-k})$.
\item $p_1 \circ q_F|_{q_F^{-1}(C)} \colon q_F^{-1}(C) \to \partial W_F$ is a $D^{n-k+1}$-bundle with structure group $\Diff(D^{n-k+1})$.
\end{enumerate}
Here, $p_1\colon C\cong \partial W_F\times [0,1]\to\partial W_F$ denotes the projection onto the first factor. 
\end{lem}

\vspace{2pt}

In particular, Lemma~\ref{lem:collar} gives a decomposition
\begin{equation}
N=q_F^{-1}(W)\cup q_F^{-1}(C), 
\end{equation}
and both pieces are disk bundles. 
Using Lemma~\ref{lem:collar} together with equation~(2.3), we obtain the following proposition.
It gives a necessary condition for Theorem\,1.4.

\vspace{2pt}

\begin{prop}\label{prop:nk1}
Let $N$ be a compact, connected $n$-dimensional manifold with boundary, and let $F\colon N\to\mathbb{R}^k$ be a boundary special generic map.
Assume that $W_F$ is contractible and that $\partial W_F$ is diffeomorphic to $S^{k-1}$.
If $n-k=1$ and $k\ge 4$, then $N\cong W_F\times D^{1}$.
\end{prop}

\proof
By Lemma\,\ref{lem:collar} and (2.3), the manifold $N$ is the union of two disk bundles $q_F^{-1}(W)$ and $q_F^{-1}(C)$.
We show that both bundles are trivial.

\vspace{2pt}

We first consider the $D^{1}$-bundle $q_F^{-1}(W)\to W$.
Since $W\cong W_F$ and $W_F$ is contractible, this bundle is trivial.
In particular,
$$
q_F^{-1}(W)\cong W_F\times D^{1}.
$$

\vspace{2pt}

We next consider the $D^{2}$-bundle $q_F^{-1}(C)\to\partial W_F$.
Since $\partial W_F\cong S^{k-1}$, its isomorphism class corresponds to an element of $\pi_{k-2}(\Diff(D^2))$. 
By~\cite{Sma}, $\Diff(D^{2})$ is homotopy equivalent to ${\rm O}(2)$.
Moreover, $H^{1}(S^{k-1};\mathbb{Z}_2)=0$ implies that $q_F^{-1}(C)$ is orientable. 
Therefore, the structure group reduces to ${\rm SO}(2)$.
Since $\pi_{k-2}({\rm SO}(2))=0$ for $k\ge 4$, the bundle is trivial.
In particular,
$$
q_F^{-1}(C)\cong S^{k-1}\times D^{2}.
$$

\vspace{2pt}

Hence, $N$ is obtained by gluing $\partial W_F\times D^{2}$ to $W_F\times D^{1}$ along ${\rm id}_{\partial W_F\times D^{1}}$.
Therefore, $N$ is diffeomorphic to $W_F\times D^{1}$.
\qed

\vspace{0.5\baselineskip}

If $k\ge 5$, then the assumptions of Proposition~\ref{prop:nk1} imply that $W_F$ is diffeomorphic to $D^k$ by \cite{Mil2}.
Hence, under the assumptions of Theorem\,1.4, $N\cong D^n$.
This proves the necessary condition in Theorem\,1.4.

\vspace{0.2\baselineskip}

We next prove the converse direction of Theorem\,1.4.
More precisely, we show the following proposition, which provides a sufficient condition for Theorem\,1.4.

\vspace{2pt}

\begin{prop}\label{prop:existence}
There exists a boundary special generic map $F\colon D^n\to\mathbb{R}^k$ whose Reeb space $W_F$ is diffeomorphic to $D^k$, where $n>k\geq 1$.
\end{prop}

\proof
Let $F$ be the composition
$$
D^{n}\hookrightarrow \mathbb{R}^{n}\xrightarrow{\pi}\mathbb{R}^k,
$$
where the first map is the standard inclusion and $\pi$ is the standard projection.
Then, $F$ is a boundary special generic map.
Moreover, the Reeb space is diffeomorphic to $D^k$. 
This follows from the definition of the first map. 
This completes the proof.
\qed

\vspace{0.5\baselineskip}

In particular, $W_F\cong D^k$ is contractible and $\partial W_F\cong S^{k-1}$. 
Hence, the sufficient condition for Theorem\,1.4 follows. 
This completes the proof of Theorem\,1.4.
\qed

\vspace{2pt}

\begin{rem}
(1)
It is not known whether a similar conclusion of Theorem\,1.4 also holds for $k=4$, that is, for $(n,k)=(5,4)$.
In the proof of the necessary condition in Theorem\,1.4, we used a result that comes from the h-cobordism theorem~\cite{Mil2}.
For this reason, the case $k=4$ must be excluded from our argument.
For a compact, contractible $4$-manifold $W$ with $\partial W\cong S^{3}$, the statement $W\cong D^{4}$ is equivalent to the smooth Poincar\'e conjecture in dimension~$4$~\cite{Mil2}.
To the best of the author's knowledge, this conjecture is still open.
If this conjecture were true, then Theorem\,1.4 would also hold for $k=4$.

\vspace{2pt}

\noindent
(2)
Proposition\,2.11 was proved under the assumptions $n-k=1$ and $k\ge 4$.
For $(n,k)=(6,4)$, we can obtain an analogous statement by a similar argument, using $\Diff(D^3)\simeq {\rm O}(3)$ by~\cite{Hat}.
More precisely, we have the following proposition.

\vspace{2pt}

\begin{prop}
Let $N$ be a compact, connected $6$-dimensional manifold with boundary, and let $F\colon N\to\mathbb{R}^4$ be a boundary special generic map.
If $W_F$ is contractible and $\partial W_F$ is diffeomorphic to $S^{3}$, then $N\cong W_F\times D^{2}$.
\end{prop}

\vspace{2pt}

However, for the same reason as in~Remark\,2.13\,(1), this does not yield a statement in the form of Theorem\,1.4.
\end{rem}

\vspace{2pt}

\begin{rem}
Lemma\,2.10 is useful for computing topological invariants of manifolds with boundary that admit boundary special generic maps.
For example, one can compute the cohomology, the Euler characteristic, and the fundamental group in terms of the Reeb space. 
In fact, they agree with those of the Reeb space.
\end{rem}

\vspace{2pt}

\section{Applications to non-singular extensions of maps}
In this section, using Theorems\,1.1--1.4 proved in Section~2, we prove Corollaries\,1.5--1.8.
We first recall the definitions of non-singular extensions and special generic maps.
We then prove Corollaries\,1.5--1.8 one by one.

\vspace{2pt}

\subsection{Definition of non-singular extensions}
We begin by recalling the notion of a non-singular extension.
Throughout this section, let $M$ a closed, connected $m$-dimensional manifold, and let $f\colon M\to\mathbb{R}^k$ a map, where $m\geq k\geq 1$.

\vspace{2pt}

\begin{defi}[\textbf{Non-singular extension}]
Assume that there exist a compact, connected ($m+1$)-dimensional manifold $N$ with $\partial N=M$, and a submersion $F\colon N\to\mathbb{R}^k$ that makes the following diagram commutative:
\begin{equation*}
  \xymatrix{
    M \ar[r]^{f} \ar[d]_{i} & \mathbb{R}^k\\
    N \ar[ru]_{F},
  }
\end{equation*}
where $i$ is the inclusion. 
Then, $F$ is called a \textit{non-singular extension} of $f$. 
\end{defi}

\vspace{2pt}

\begin{rem}
In this paper, a non-singular extension means a submersion on a manifold with boundary whose restriction to the boundary coincides with the given map.
However, in~\cite{Cur, Iwa, Sei}, a non-singular extension is defined as a submersion obtained by extending a given submersion defined on a collar neighborhood of the boundary.
This difference is only technical.
Indeed, we typically first extend the given map to a submersion on a collar, and then extend it further over the whole manifold; the previous works incorporate the collar step into the definition.
\end{rem}

\vspace{2pt}

\subsection{Definition of special generic maps}
We recall the definition of special generic maps on closed manifolds.
For their basic properties, see \cite{Sae1, SS1, SS2}. 
We first recall the notion of a definite fold point.

\vspace{2pt}

\begin{defi}[\textbf{Definite fold point}]
A point $p\in M$ is called a \textit{definite fold point} of $f$ if there exist local coordinates $(x_1,\dots,x_m)$ of $M$ around $p$ and $(y_1,\dots,y_k)$ of $\mathbb{R}^k$ around $f(p)$ such that $f$ is written as 
$$
f; (x_1,\dots,x_m)\mapsto (y_1,\dots,y_k)=\Big(x_1,\dots,x_{k-1},\sum_{i=k}^{m}x_i^2\Big). 
$$
\end{defi}

\begin{defi}[\textbf{Special generic map}]
The map $f$ is called a \textit{special generic map} if all singular points of $f$ are definite fold points. 
\end{defi}

\vspace{2pt}

\begin{rem}
By Remark\,2.5 and Definition\,3.1, every boundary special generic map gives a non-singular extension of a special generic map.
On the other hand, the converse does not hold in general. 
Namely, even if a special generic map admits a non-singular extension, such an extension need not be a boundary special generic map.
We present such an example in Remark\,3.9 and Figure\,1.
\end{rem}

\vspace{2pt}

\subsection{Proof of Corollary\,1.5}
In this subsection, we prove Corollary\,1.5 by applying Theorem\,1.1. 
For a special generic map into $\mathbb{R}$, the following statements are equivalent: admitting a boundary special generic map as its non-singular extension, and merely admitting a non-singular extension.
Indeed, if a special generic map into $\mathbb{R}$ admits a non-singular extension, then the collar neighborhoods of the boundary around the local maximum and local minimum satisfy the condition of a boundary definite fold point.
This follows from the compactness of the source manifold together with the description of the fiber changes near these extrema.
Therefore, it suffices to show that the following are equivalent: a special generic map $f$ into $\mathbb{R}$ admits a boundary special generic map as a non-singular extension, and $M \cong S^m$.

\vspace{0.2\baselineskip}

Assume that $f$ admits a boundary special generic map $F\colon N\to\mathbb{R}$ as a non-singular extension.
Then, $N$ is diffeomorphic to $D^{m+1}$ by Theorem\,1.1. 
Hence, $M=\partial N\cong S^m$.

\vspace{0.2\baselineskip}

Conversely, assume that $M\cong S^m$.
By~\cite[Theorem 4.1]{ST}, $f$ admits an immersion lift to $\mathbb{R}^{m+1}$; that is, there exists an immersion $\eta\colon M\to\mathbb{R}^{m+1}$ such that $f=\pi\circ\eta$, where $\pi\colon\mathbb{R}^{m+1}\to\mathbb{R}$ is the standard projection.
Taking a collar neighborhood of $M$, we obtain a submersion $\tilde f\colon S^m\times[0,1)\to\mathbb{R}$ such that, at the local maximum point of $f$, the image under $d\tilde f$ of an outward normal vector points in the positive direction of $\mathbb{R}$, while at the local minimum point it points in the negative direction.
Then, $\tilde f$ extends to a submersion $D^{m+1}\to\mathbb{R}$ by~\cite[Proposition 2.3]{Sei}.
By the equivalence stated above, this completes the proof of Corollary\,1.5.
\qed

\vspace{0.5\baselineskip}

\begin{rem}
For maps of closed manifolds, the existence of immersion and embedding lifts has been studied as a fundamental problem~\cite{Hae, Nis, ST} of singularity theory.
We mean an embedding lift of $f$ to $\mathbb{R}^{m+1}$ as an embedding $\eta\colon M\to\mathbb{R}^{m+1}$ such that $f=\pi\circ\eta$, where $\pi\colon\mathbb{R}^{m+1}\to\mathbb{R}$ is the standard projection.
For a special generic map $f$ into $\mathbb{R}$, Saeki--Takase~\cite{ST} show that it admits an embedding lift to $\mathbb{R}^{m+1}$ if and only if $M\cong S^m$.
Therefore, together with Corollary\,1.5, we obtain that, for special generic maps into $\mathbb{R}$, the existence of a non-singular extension is equivalent to the existence of an embedding lift to $\mathbb{R}^{m+1}$.
\end{rem}

\subsection{Proof of Corollary\,1.6}
In this subsection, we prove Corollary\,1.6. 
We first record the following lemma, which follows directly from Theorem\,1.2.

\vspace{2pt}

\begin{lem}
Let $f\colon M\to\mathbb{R}^2$ be a special generic map. 
If $f$ admits a boundary special generic map as a non-singular extension, then there exist integers $r,r'\geq 0$ such that 
$$
M\cong(\sharp^r S^1\times S^{m-1})\,\sharp\,(\sharp^{r'} S^1\tilde{\times} S^{m-1}), 
$$
where $S^1\tilde{\times} S^{m-1}$ denotes the non-trivial $S^{m-1}$-bundle over $S^1$. 
Here, we interpret the case of zero summands as $S^{m}$.
\end{lem}

\vspace{2pt}

We first show that a special generic map $f\colon \Sigma^m\to\mathbb{R}^2$ does not admit a boundary special generic map as a non-singular extension, where $\Sigma^m$ is an exotic sphere. 
Assume that $f$ admits such a non-singular extension.
Then, by Lemma\,3.7 together with the discussion on fundamental groups, we obtain that $\Sigma^m$ is diffeomorphic to the standard sphere $S^m$.
This contradicts the definition of an exotic sphere.

\vspace{2pt}

We next prove that a special generic map $f\colon (S^1\times S^{m-1})\,\sharp\,\Sigma^m\to\mathbb{R}^2$ does not admit a boundary special generic map as a non-singular extension. 
Assume that $f$ admits such a non-singular extension. 
By Lemma\,3.7, the source manifold is diffeomorphic to $(\sharp^r S^1\times S^{m-1})\,\sharp\,(\sharp^{r'} S^1\tilde{\times} S^{m-1})$ for some $r,r'\geq 0$. 
By considering orientability and the fundamental group, we obtain $r=1$ and $r'=0$. 
Hence, 
$$
(S^1\times S^{m-1})\,\sharp\,\Sigma^m\cong S^1\times S^{m-1}. 
$$
On the other hand, by \cite{Sch}, the inertia group $I(S^1\times S^{m-1})$ is trivial.
In particular, when $S^1\times S^{m-1}$ is oriented, there is no oriented exotic sphere $\Sigma^m$ such that $(S^1\times S^{m-1})\,\sharp\,\Sigma^m$ is orientation-preserving diffeomorphic to $S^1\times S^{m-1}$.
Moreover, $S^1\times S^{m-1}$ is orientation-preserving diffeomorphic to $\overline{S^1\times S^{m-1}}$.
Therefore, the above diffeomorphism cannot exist, and the assumption leads to a contradiction.
This completes the proof of Corollary\,1.6. 
\qed

\vspace{0.5\baselineskip}

\begin{rem}
(1) 
Saeki~\cite{Sae1} completely classified the diffeomorphism types of closed manifolds that admit special generic maps into $\mathbb{R}^2$.
Besides the manifolds considered in Corollary\,1.6, Saeki's list also contains $S^1\tilde{\times}S^{m-1}$, $S^1\times\Sigma^{m-1}$, $S^1\tilde{\times}\Sigma^{m-1}$, and their connected sums.
It is unknown whether special generic maps on these manifolds into $\mathbb{R}^2$ admit a boundary special generic map as a non-singular extension.

\vspace{2pt}

\noindent
(2) On the other hand, the inertia group of $S^1\times \Sigma^{m-1}$ was computed by Kawakubo~\cite{Kaw}. 
In particular, when this group is trivial for some $m$, one can show the non-existence of such extensions by applying the same argument as in the proof above. 
\end{rem}

\begin{rem}
In the proof of Corollary\,1.6, unlike in the proof of Corollary\,1.5, we cannot use the equivalence between the following two conditions: that a special generic map admits a boundary special generic map as a non-singular extension, and that it admits a non-singular extension.
In fact, as shown in Figure\,1, there exists a special generic map $S^2\to\mathbb{R}^2$ that admits a non-singular extension but does not admit a boundary special generic map as a non-singular extension.
We also note that it is still an open question whether every special generic map into $\mathbb{R}^2$ admits a non-singular extension.

\begin{figure}[t]
\centering
\includegraphics[width=40mm]{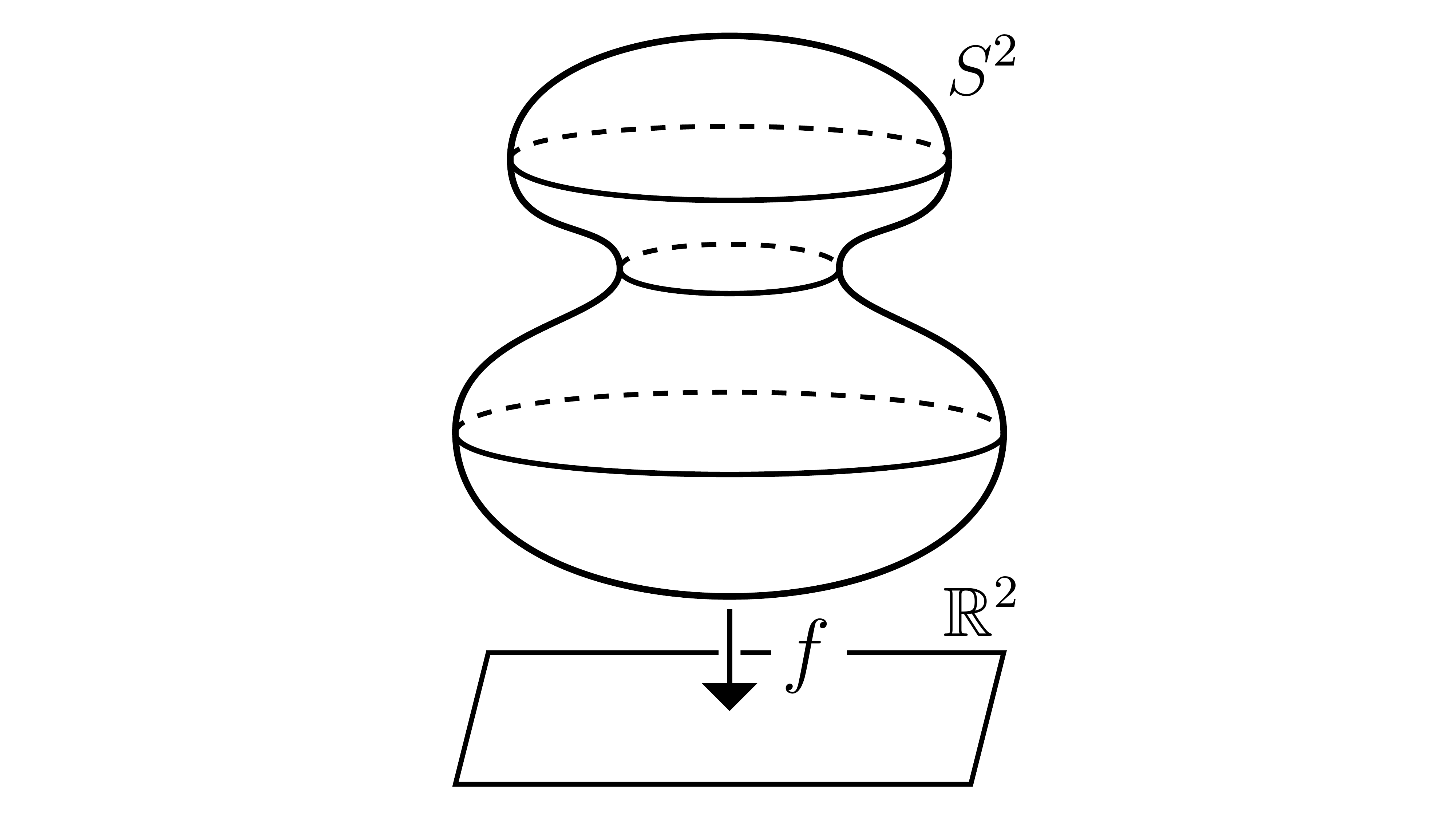}
\caption{
A special generic map $f\colon S^2\to\mathbb{R}^2$ that admits a non-singular extension but admits no boundary special generic map as a non-singular extension.
}
\end{figure}

\end{rem}

\vspace{2pt}

\subsection{Proof of Corollary\,1.7}
In this subsection, we prove Corollary\,1.7. 
Before giving the proof, we record the following lemma, which follows immediately from Theorem\,1.3.

\vspace{2pt}

\begin{lem}
Let $f\colon M\to\mathbb{R}^3$ be a special generic map of a $3$-manifold $M$. 
If $f$ admits a boundary special generic map as a non-singular extension whose source manifold is simply connected, then there exist integers $r\geq 0$ and $k_1,\dots,k_r$ such that 
$$
M\cong L(k_1,1)\,\sharp\dots\sharp\,L(k_r,1), 
$$
where $L(k,1)$ is a lens space. 
Here, we interpret the case of zero summands as $S^3$.
\end{lem}

\vspace{2pt}

Let $M$ be a closed, connected, irreducible $3$-manifold with a non-cyclic fundamental group, and let $f\colon M\to\mathbb{R}^3$ be a special generic map.
Assume that $f$ admits a boundary special generic map as a non-singular extension whose source manifold is simply connected. 
By Lemma\,3.10, there exist integers $r\ge 0$ and $k_1,\dots,k_r$ such that
$$
M\cong L(k_1,1)\,\sharp\dots\sharp\,L(k_r,1). 
$$
If $r=0$, then $M\cong S^3$. 
However, it contradicts the assumption that $\pi_1(M)$ is non-cyclic.
Hence, $r\ge 1$.
Since $M$ is irreducible, it follows that $M\cong L(k_i,1)$ for some $i\in\{1,\dots,r\}$.
However, $\pi_1(M)$ is non-cyclic, whereas $\pi_1(L(k_i,1))$ is cyclic.
This is a contradiction. 
\qed

\vspace{0.5\baselineskip}

\begin{rem}
(1) 
Corollary\,1.7 assumes that $M$ is irreducible and has a non-cyclic fundamental group only to rule out the case where $M$ is a connected sum of lens spaces.
Accordingly, we may replace these assumptions by any condition implying that $M$ is not diffeomorphic to a connected sum of lens spaces.
For example, it is enough to assume that $H_2(M;\mathbb{Z})\neq 0$.

\vspace{2pt}

\noindent
(2)
In Corollary\,1.7, we used Lemma\,3.10, which follows from the $3$-dimensional case of Theorem\,1.3. 
For $4$-manifolds, one also obtains the following analogue of Lemma\,3.10.

\vspace{2pt}

\begin{prop} 
Let $f\colon M\to\mathbb{R}^3$ be a special generic map of a $4$-manifold $M$. 
If $f$ admits a boundary special generic map as a non-singular extension whose source manifold is simply connected, then there exists an integer $k\geq 0$ such that 
$$
M\cong \sharp^k(\mathbb{CP}^2 \,\sharp\,\overline{\mathbb{CP}^2})
\text{\,\,\,or\,\,\,}
\sharp^k(S^2 \times S^2). 
$$
Here, we interpret the case of zero summands as $S^4$.
\end{prop}

\vspace{2pt}

This follows from the fact by~\cite{Kim} that $M$ is diffeomorphic to the double of a $4$-dimensional $2$-handlebody, together with the classification of such doubles in~\cite{GS}.
Moreover, which of the two types occurs is determined by whether the intersection form $Q_M$ is even or not.

\vspace{2pt}

\noindent
(3) 
In higher dimensions, one can obtain similar consequences by combining Theorem\,1.3 with the result in~\cite{Kim} as follows.

\vspace{2pt}

\begin{prop}
Let $f\colon M\to\mathbb{R}^3$ be a special generic map of an $m$-dimensional manifold $M$, where $m\ge 7$, or $m=5,6$ with $\sharp S(f)\leq 2$.
If $f$ admits a boundary special generic map as a non-singular extension whose source manifold is simply connected, then $M$ is an $m$-dimensional $2$-handlebody.
\end{prop}
\end{rem}

\vspace{2pt}

We next state a consequence of Corollary\,1.7 for non-singular extensions of special generic maps on Brieskorn manifolds.
$$
\Sigma(a,b,c)=\{\,(x,y,z)\in\mathbb{C}^3\mid x^a+y^b+z^c=0\,\}\cap S^5
$$ 
is called a $3$-dimensional \textit{Brieskorn manifold}, where $a,b,c\in\mathbb{Z}_{\geq 2}$. 
Since $\Sigma(a,b,c)$ is the link of the isolated hypersurface singularity, it is irreducible~\cite{Neu}.
It is also known that, for many triples $(a,b,c)$, the group $\pi_1(\Sigma(a,b,c))$ is non-cyclic~\cite{Mil3}.
For example, this holds for $(a,b,c)=(2,5,7),(3,4,5),(2,3,13)$.
We denote the corresponding Brieskorn manifolds by $\Sigma_1,\Sigma_2,\Sigma_3$.
Moreover, each $\Sigma_i$ is known to arise as the boundary of a compact, contractible $4$-manifold~\cite{AK}.
We denote such a manifold by $W_i$, so that $\partial W_i=\Sigma_i$.
Then, Corollary\,1.7 immediately yields the following corollary.

\vspace{2pt}

\begin{cor}
Assume that a special generic map from $\Sigma_i$ to $\mathbb{R}^3$ admits a non-singular extension from $W_i$ to $\mathbb{R}^3$. 
Then, no such non-singular extension can be a boundary special generic map. 
\end{cor}

\vspace{2pt}

\begin{rem}
Casson--Harer~\cite{CH} also constructed many Mazur-type manifolds whose boundaries are Brieskorn manifolds.
Since the argument above uses only that $\pi_1(\Sigma(a,b,c))$ is non-cyclic, the same argument applies to the examples in~\cite{CH} whenever $\pi_1(\Sigma(a,b,c))$ is non-cyclic.
\end{rem}

\vspace{2pt}

\subsection{Proof of Corollary\,1.8}
In this subsection, we prove Corollary\,1.8.
Let $f\colon \Sigma^m\to\mathbb{R}^m$ be a special generic map, and assume that $f$ admits a boundary special generic map as a non-singular extension whose Reeb space is contractible.
Then, Theorem\,1.4 implies that $\Sigma^m$ is diffeomorphic to $S^m$.
However, since $\Sigma^m$ is an exotic sphere, this contradicts the definition.
This completes the proof of Corollary\,1.8.
\qed

\vspace{0.5\baselineskip}

\begin{rem}
Corollary\,1.8 shows that $\Sigma^m$ admits no boundary special generic map with contractible Reeb space as a non-singular extension.
On the other hand, without the assumption on the Reeb space, it is still unknown whether a special generic map $\Sigma^m\to\mathbb{R}^m$ admits a boundary special generic map as a non-singular extension.
\end{rem}

\vspace{0.2\baselineskip}

\section*{Aknowledgement}
The author would like to thank Osamu Saeki, Noriyuki Hamada, and Naoki Kitazawa for their useful discussions and comments. 
This work has been partially supported by JSPS KAKENHI Grant Number JP23H05437 and by WISE program (MEXT) at Kyushu University.

\end{document}